\documentclass[12pt,a4paper]{article}
\usepackage[dvips]{graphicx}
\usepackage{makeidx}
\usepackage{amsmath}
\usepackage{amsfonts}

\newcommand {\ess}{{\mathrm{ess}}}

\newtheorem{theorem}{Theorem}

\newtheorem{assume}{Assumption}

\newtheorem{hyp}{H}
\newtheorem{problem}{Problem}

\newcommand {\epsvec}{{\boldsymbol{\epsilon}}}

\newcommand {\omegavec}{{\boldsymbol{\omega}}}

\newcommand {\alphavec}{\boldsymbol{\alpha}}
\newcommand {\xivec}{{\boldsymbol{\xi}}}

\newcommand {\phivec}{\mbox{\boldmath $\phi$}}

\newcommand {\thetavec}{{\boldsymbol{\theta}}}

\newcommand{\bfx}{\mathbf{x}}
\newcommand{\bff}{\mathbf{f}}
\newcommand{\bfz}{\mathbf{z}}

\newcommand{\bfy}{\mathbf{y}}
\newcommand{\bfh}{\mathbf{h}}

\newcommand{\bfu}{\mathbf{u}}
\newcommand{\bfg}{\mathbf{g}}

\newcommand{\bfq}{\mathbf{q}}

\newcommand{\dpsi}{{\dot\psi}}
\newcommand{\pd}{{\partial}}

\newcommand{\Real}{\mathbb{R}}

\begin{document}

\baselineskip 0.6cm

\title{\bf Decentralized adaptation in interconnected uncertain systems with nonlinear parametrization}
\author{Ivan Tyukin\thanks{({\bf corresponding author}) Laboratory for Perceptual Dynamics, Brain Science Institute, RIKEN (Institute for Physical and Chemical Research)
                           , 2-1, Hirosawa, Wako-shi, Saitama, 351-0198, Japan, e-mail:
                           \{tyukinivan\}@brain.riken.jp},  Cees van Leeuwen\thanks{Laboratory for Perceptual Dynamics, Brain
Science Institute, RIKEN (Institute for Physical and Chemical
Research)
                          , 2-1, Hirosawa, Wako-shi, Saitama, 351-0198, Japan, e-mail:
                           \{ceesvl\}@brain.riken.jp}
}

\date{}

\maketitle{}

\begin{abstract}
We propose a technique for the design and analysis of
decentralized adaptation algorithms in interconnected dynamical
systems. Our technique does not require Lyapunov stability of the
target dynamics and allows nonlinearly parameterized
uncertainties. We show that for the considered class of systems,
conditions for reaching the control goals can be formulated in
terms of the nonlinear $L_2$-gains of target dynamics of each
interconnected subsystem. Equations for decentralized controllers
and corresponding adaptation algorithms are also explicitly
provided.

{\it Keywords:} nonlinear parametrization; unstable,
non-equilibrium dynamics; decentralized adaptive control; monotone
functions
\end{abstract}

\section*{Notation}

According to the standard convention,  $\Real$ defines the field
of real numbers and  $\Real_{\geq c}=\{x\in\Real|x\geq c\}$,
$\Real_{+}=\Real_{\geq 0}$; symbol $\Real^n$ stands for a linear
space $\mathcal{L}(\Real)$ over the field of reals with
$\mathrm{dim}\{\mathcal{L}(\Real)\}=n$; $\|\bfx\|$ denotes the
Euclidian norm of $\bfx\in\Real^n$; $\mathcal{C}^k$ denotes the
space of functions that are  at least $k$ times differentiable;
$\mathcal{K}$ denotes the class of all strictly increasing
functions $\kappa: \Real_+\rightarrow \Real_+$ such that
$\kappa(0)=0$. By  ${L}_{p}^n[t_0,T]$, where $T>0$, $p\geq 1$ we
denote the space of all functions $\bff:\Real_+\rightarrow\Real^n$
such that
$\|\bff\|_{p,[t_0,T]}=\left(\int_{0}^T\|\bff(\tau)\|^{p}d\tau\right)^{1/p}<\infty$;
$\|\bff\|_{p,[t_0,T]}$ denotes the ${L}_{p}^n[t_0,T]$-norm of
$\bff(t)$. By ${L}^n_\infty[t_0,T]$ we denote the space of all
functions $\bff:\Real_+\rightarrow\Real^n$ such that
$\|\bff\|_{\infty,[t_0,T]}=\ess \sup\{\|\bff(t)\|,t \in
[t_0,T]\}<\infty$, and $\|\bff\|_{\infty,[t_0,T]}$ stands for the
${L}^n_\infty[t_0,T]$ norm of $\bff(t)$.

A function $\bff(\bfx): \Real^{n}\rightarrow \Real^m$ is said to
be locally bounded if for any $\|\bfx\|<\delta$ there exists a
constant $D(\delta)>0$ such that the following inequality holds:
$\|\bff(\bfx)\|\leq D(\delta)$. Let $\Gamma$ be an $n\times n$
square matrix, then $\Gamma>0$ denotes a positive definite
(symmetric) matrix, and $\Gamma^{-1}$ is the inverse of $\Gamma$.
By $\Gamma\geq 0$ we denote a positive semi-definite matrix,
$\|\bfx\|_{\Gamma}^2$ to denotes the quadratic form:
$\bfx^{T}\Gamma\bfx$, $\bfx\in\Real^n$. The notation $|\cdot|$
stands for the modulus of a scalar. The solution of a system of
differential equations $\dot{\bfx}=\bff(\bfx,t,\thetavec,\bfu), \
\bfx(t_0)=\bfx_0$, $\bfu:\Real_+\rightarrow\Real^m$,
$\thetavec\in\Real^d$ for $t\geq t_0$ will be denoted as
$\bfx(t,\bfx_0,t_0,\thetavec,\bfu)$, or simply as $\bfx(t)$ if it
is clear from the context what  the values of $\bfx_0,\thetavec$
are and how the function $\bfu(t)$ is defined.

Let $\bfu:\Real^n\times\Real^d\times\Real_+\rightarrow\Real^m$ be a
function of state $\bfx$, parameters $\hat{\thetavec}$, and time
$t$. Let in addition both $\bfx$ and $\hat{\thetavec}$ be functions
of $t$. Then in case the arguments of $\bfu$ are clearly defined by
the context,  we will simply write $\bfu(t)$ instead of
$\bfu(\bfx(t),\hat{\thetavec}(t),t)$.

The (forward complete) system
$\dot{\bfx}=\bff(\bfx,t,\thetavec,\bfu(t))$, is said to have an
$L_{p}^m [t_0,T]\mapsto L_{q}^n[t_0,T]$, gain ($T\geq t_0$,
$p,q\in\Real_{\geq 1}\cup\infty$) with respect to its input
$\bfu(t)$ if and only if $\bfx(t,\bfx_0,t_0,\thetavec,\bfu(t))\in
L_{q}^n [t_0,T]$ for any $\bfu(t)\in L_{p}^m [t_0,T]$ and there
exists a function
$\gamma_{q,p}:\Real^n\times\Real^d\times\Real_+\rightarrow\Real_+$
such that the following inequality holds:
$\|\bfx(t)\|_{q,[t_0,T]}\leq
\gamma_{q,p}(\bfx_0,\thetavec,\|\bfu(t)\|_{p,[t_0,T]})$. The
function $\gamma_{q,p}(\bfx_0,\thetavec,\|\bfu(t)\|_{p,[t_0,T]})$
is assumed to be non-decreasing  in $\|\bfu(t)\|_{p,[t_0,T]}$, and
locally bounded in its arguments.

For notational convenience when dealing with vector fields and
partial derivatives we will use the following extended notion of
the Lie derivative of a function. Let $\bfx\in\Real^n$ and assume
$\bfx$ can be partitioned as follows $\bfx=\bfx_1\oplus\bfx_2$,
where $\bfx_1\in\Real^q$, $\bfx_1=(x_{11},\dots,x_{1q})^T$,
$\bfx_2\in\Real^p$, $\bfx_2=(x_{21},\dots,x_{2p})^T$, $q+p=n$, and
$\oplus$ denotes the  concatenation of two vectors. Define
$\bff:\Real^{n}\rightarrow\Real^n$ such that
$\bff(\bfx)=\bff_1(\bfx)\oplus\bff_2(\bfx)$, where
$\bff_1:\Real^n\rightarrow\Real^q$,
$\bff_1(\cdot)=(f_{11}(\cdot),\dots,f_{1q}(\cdot))^T$,
$\bff_2:\Real^n\rightarrow\Real^p$,
$\bff_2(\cdot)=(f_{21}(\cdot),\dots,f_{2p}(\cdot))^T$. Then
$L_{\bff_i(\bfx)}\psi(\bfx,t)$, $i\in\{1,2\}$ denotes the Lie
derivative of the function $\psi(\bfx,t)$ with respect to the
vector field $\bff_i(\bfx,\thetavec)$:
$L_{\bff_i(\bfx)}\psi(\bfx,t)=\sum_{j}^{\dim{\bfx_i}}\frac{\pd
\psi(\bfx,t) }{\pd x_{ij}}f_{ij}(\bfx,\thetavec)$.

\section{Introduction}

We consider the problem how to control the behavior of complex
dynamical systems composed of interconnected lower-dimensional
subsystems. Centralized control of these systems is practically
inefficient because of high demands for computational power,
measurements  and prohibitive communication cost. On the other
hand, standard decentralized solutions often face severe
limitations due to the deficiency of information about the
interconnected subsystems. In addition, the nature of their their
interconnections may vary depending on conditions in the
environment. In order to address these problems in their most
general setup, decentralized adaptive control is needed.

Currently there is a large literature on decentralized adaptive
control which contains successful solutions to problems of
adaptive stabilization \cite{Gavel_1989,Jain_1997}, tracking
\cite{Ioannou86,Jain_1997,Shi_1992,Passino96}, and output
regulation \cite{Jiang_2000,Huang_2003} of linear and nonlinear
systems. In most of these cases the problem of decentralized
control is solved within the conventional framework of adaptive
stabilization/tracking/regulation by a family of linearly
parameterized controllers. While these results may be successfully
implemented in a large variety of technical and artificial
systems, there is room for further improvements. In particular,
when the target dynamics of the systems is not stable in the
Lyapunov sense but intermittent, meta-stable, or multi-stable
\cite{Arecchi_2004,Raffone_2003,Tsuda_2004} or when the
uncertainties are nonlinearly parameterized
\cite{Armstrong_1993,Boskovic_1995,Canudas_1999,Kitching_2000},
and no domination of the uncertainties by feedback is allowed.

In the present article we address these issues at once for a class
of nonlinear dynamical systems. Our contribution is that we
provide conditions ensuring forward-completeness, boundedness and
asymptotic reaching of the goal for a pair of interconnected
systems with uncertain coupling and parameters. Our method does
not require availability of a Lyapunov function for the desired
motions in each subsystem, nor linear parametrization of the
controllers. Our results can straightforwardly be extended to
interconnection of arbitrary many (but still, a finite number of)
subsystems. Explicit equations for corresponding decentralized
adaptive controllers are also provided.

The paper is organized as follows. In Section 2 we provide a
formal statement of the problem, Section 3 contains necessary
preliminaries and auxiliary results. In Section 4 we present the
main results of our current contribution,  and in Section 5 we
provide concluding remarks to our approach.

\section{Problem Formulation}

Let us consider two interconnected systems $\mathcal{S}_x$ and
$\mathcal{S}_y$:
\begin{eqnarray}
&\mathcal{S}_x: & \
\dot{\bfx}=\bff(\bfx,\thetavec_x)+\gamma_y(\bfy,t)+
\bfg(\bfx)u_x \label{eq:system:s1} \\
&\mathcal{S}_y: & \
\dot{\bfy}=\bfq(\bfy,\thetavec_y)+\gamma_x(\bfx,t)+\bfz(\bfy)u_y\label{eq:system:s2}
\end{eqnarray}
where $\bfx\in\Real^{n_x}$, $\bfy\in\Real^{n_y}$ are the state
vectors of systems $\mathcal{S}_x$ and $\mathcal{S}_y$, vectors
$\thetavec_x\in\Real^{n_{\theta_x}}$,
$\thetavec_y\in\Real^{n_{\theta_y}}$ are unknown parameters,
functions
$\bff:\Real^{n_x}\times\Real^{n_{\theta_x}}\rightarrow\Real^{n_x}$,
$\bfq:\Real^{n_y}\times\Real^{n_{\theta_y}}\rightarrow\Real^{n_y}$,
$\bfg:\Real^{n_x}\rightarrow\Real^{n_x}$,
$\bfz:\Real^{n_y}\rightarrow\Real^{n_y}$ are continuous and
locally bounded. Functions
$\gamma_y:\Real^{n_y}\times\Real_+\rightarrow\Real_n$,
$\gamma_x:\Real^{n_x}\times\Real_+\rightarrow\Real^{n_y}$, stand
for nonlinear, non-stationary and, in general, unknown couplings
between systems $\mathcal{S}_x$ and $\mathcal{S}_y$, and
$u_x\in\Real$, $u_y\in\Real$ are the control inputs.

In the present paper we are interested in the following problem

\begin{problem}\label{problem:decentralized}\normalfont Let $\psi_x:\Real^{n_x}\times\Real_+\rightarrow\Real$,
$\psi_y:\Real^{n_y}\times\Real_+\rightarrow\Real$ be the goal
functions for systems $\mathcal{S}_x$, $\mathcal{S}_y$
respectively. In the other words, for some values
$\varepsilon_x\in\Real_{+}$, $\varepsilon_y\in\Real_+$ and time
instant $t^\ast\in\Real_+$, inequalities
\begin{equation}\label{eq:goal_functionals}
 \|\psi_x(\bfx(t),t)\|_{\infty,[t^\ast,\infty]}\leq\varepsilon_x, \
 \|\psi_y(\bfy(t),t)\|_{\infty,[t^\ast,\infty]}\leq\varepsilon_y
\end{equation}
specify the desired state of interconnection (\ref{eq:system:s1}),
(\ref{eq:system:s2}). Derive functions $u_x(\bfx,t)$, $u_y(\bfy,t)$
such that for all  $\thetavec_x\in\Real^{n_{\theta_x}}$,
$\thetavec_y\in\Real^{n_{\theta_y}}$

1) interconnection (\ref{eq:system:s1}), (\ref{eq:system:s2}) is
forward-complete;

2) the trajectories $\bfx(t)$, $\bfy(t)$ are bounded;

3) for given values of $\varepsilon_x$, $\varepsilon_y$, some
$t^\ast\in\Real_+$ exists such that inequalities
(\ref{eq:goal_functionals}) are satisfied or, possibly, both
functions $\psi_x(\bfx(t),t)$, $\psi_y(\bfy(t),t)$ converge to
zero as $t\rightarrow\infty$.

Function $u_x(\cdot)$ should not depend explicitly on $\bfy$ and,
symmetrically, function $u_y(\cdot)$ should not depend explicitly
on $\bfx$. The general structure of the desired configuration of
the control scheme is provided in Figure 1.
\end{problem}

\begin{figure}
\begin{center}
\includegraphics[width=110pt]{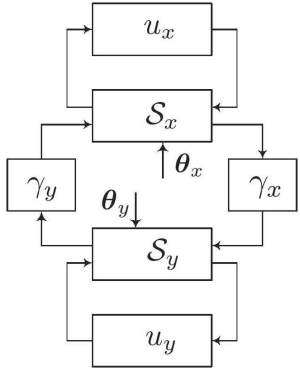}
\end{center}
\begin{center}
\caption{General structure of
interconnection}\label{fig:decentralized:singularity}
\end{center}
\end{figure}

In the next sections we provide sufficient conditions, ensuring
solvability of Problem \ref{problem:decentralized} and we also
explicitly derive functions $u_x(\bfx,t)$ and $u_y(\bfy,t)$ which
satisfy requirements 1) -- 3) of Problem
\ref{problem:decentralized}. We start with the introduction of a
new class of adaptive control schemes and continue by providing
the input-output characterizations of the controlled systems.
These results are given in Section \ref{sec:preliminary}. Then,
using these characterizations, in Section \ref{sec:main} we
provide the main results of our study.

\section{Assumptions and properties of the decoupled systems}\label{sec:preliminary}

Let the following system be given:
\begin{equation}\label{system1}
\begin{split}
\dot{\bfx}_1=&\bff_1(\bfx)+\bfg_1(\bfx)u, \\
\dot{\bfx}_2=&\bff_2(\bfx,\thetavec)+\bfg_2(\bfx)u,
\end{split}
\end{equation}
where
\[
\bfx_1=(x_{11},\dots,x_{1 q})^T\in \Real^q; \
\bfx_2=(x_{21},\dots,x_{2 p})^T\in \Real^p;
\]
\[
\bfx=(x_{11},\dots,x_{1 q},x_{21},\dots,x_{2 p})^T\in \Real^{n}
\]
$\thetavec\in \Omega_\theta\in \Real^d$ is a vector of unknown
parameters, and $\Omega_\theta$ is a closed bounded subset of
$\Real^d$; $u\in\Real$ is the control input, and functions
$\bff_1:\Real^{n}\rightarrow \Real^{q}$,
$\bff_2:\Real^{n}\times\Real^d\rightarrow \Real^{p}$,
$\bfg_1:\Real^{n}\rightarrow \Real^q$,
$\bfg_2:\Real^{n}\rightarrow\Real^{p}$ are continuous and locally
bounded. The vector $\bfx\in\Real^n$ is the state vector, and
vectors $\bfx_1$, $\bfx_2$ are referred to as {\it
uncertainty-independent} and {\it uncertainty-dependent} partition
of $\bfx$, respectively. For the sake of compactness we will also
use the following description of (\ref{system1}):
\begin{equation}\label{system}
\dot{\bfx}=\bff(\bfx,\thetavec)+\bfg(\bfx)u,
\end{equation}
where
\[
\bfg(\bfx)=(g_{11}(\bfx),\dots,g_{1q}(\bfx),g_{21}(\bfx),\dots,g_{2
p}(\bfx))^{T},
\]
\[
\bff(\bfx)=(f_{11}(\bfx),\dots,f_{1q}(\bfx),f_{21}(\bfx,\thetavec),\dots,f_{2
p}(\bfx,\thetavec))^{T}.
\]

As a measure of closeness of trajectories $\bfx(t)$ to the desired
state we introduce the error or goal function $\psi:\Real^n\times
\Real_+\rightarrow \Real, \ \psi\in \mathcal{C}^1$.
We suppose also that for the chosen function $\psi(\bfx,t)$
satisfies the following:
\begin{assume}[Target operator]\label{assume:psi}  For the given function $\psi(\bfx,t)\in \mathcal{C}^1$ the following  property holds:
\begin{equation}\label{eq:assume_psi}
\|\bfx(t)\|_{\infty,[t_0,T]}\leq
\tilde{\gamma}\left(\bfx_0,\thetavec,\|\psi(\bfx(t),t)\|_{\infty,[t_0,T]}\right)
\end{equation}
where
$\tilde{\gamma}\left(\bfx_0,\thetavec,\|\psi(\bfx(t),t)\|_{\infty,[t_0,T]}\right)$
is a locally bounded and non-negative function of its arguments.
\end{assume}
Assumption \ref{assume:psi} can be interpreted as a sort of {\it
unboundedness observability} property \cite{Jiang_1994} of system
(\ref{system1}) with respect to the ``output" function
$\psi(\bfx,t)$. It can also be viewed as a {\it bounded input -
bounded state} assumption for system (\ref{system1}) along the
constraint
$\psi(\bfx(t,\bfx_0,t_0,\thetavec,u(\bfx(t),t)),t)=\upsilon(t)$,
where the signal $\upsilon(t)$ serves as a new input. If, however,
boundedness of the state is not explicitly required  (i.e. it is
guaranteed by additional control or follows from the physical
properties of the system itself), Assumption \ref{assume:psi} can
be removed from the statements of our results.

Let us specify a class of control inputs $u$ which can ensure
boundedness of $\bfx(t,\bfx_0,t_0,\thetavec,u)$ for every
$\thetavec\in \Omega_\theta$ and $\bfx_0\in\Real^n$. According to
(\ref{eq:assume_psi}), boundedness of
$\bfx(t,\bfx_0,t_0,\thetavec,u)$ is ensured if we find a control
input $u$ such that $\psi(\bfx(t),t)\in L_\infty^1[t_0,\infty]$.
For this objective consider the dynamics of system (\ref{system})
with respect to $\psi(\bfx,t)$:
\begin{equation}\label{dpsi}
\dot{\psi}=L_{\bff(\bfx,\thetavec)}\psi(\bfx,t)+L_{\bfg(\bfx)}\psi(\bfx,t)u+\frac{\pd
\psi(\bfx,t)}{\pd t},
\end{equation}
Assuming that the inverse
$\left(L_{\bfg(\bfx)}\psi(\bfx,t)\right)^{-1}$ exists everywhere,
we may choose the control input $u$ in the following class of
functions:
\begin{equation}\label{control}
\begin{split}
u(\bfx,\hat\thetavec,\omegavec,t)&=\frac{1}{L_{\bfg(\bfx)}\psi(\bfx,t)}\left(-L_{\bff(\bfx,\hat{\thetavec})}\psi(\bfx,t)-\varphi(\psi,\omegavec,t)-\frac{\pd\psi(\bfx,t)}{\pd
t}\right) \\
& \ \varphi: \ \Real\times\Real^w\times\Real_+\rightarrow\Real
\end{split}
\end{equation}
where $\omegavec\in\Omega_\omega\subset\Real^w$ is a vector of
{\it known} parameters of the function
$\varphi(\psi,\omegavec,t)$. Denoting
$L_{\bff(\bfx,\thetavec)}\psi(\bfx,t)=f(\bfx,\thetavec,t)$ and
taking into account (\ref{control}) we may rewrite equation
(\ref{dpsi}) in the following manner:
\begin{equation}\label{error_model}
\dpsi=f(\bfx,\thetavec,t)-f(\bfx,\hat{\thetavec},t)-\varphi(\psi,\omegavec,t)
\end{equation}

For the purpose of the present article, instead of
(\ref{error_model}) it is worthwhile  to consider the extended
equation:
\begin{equation}\label{error_model_d}
\dpsi=f(\bfx,\thetavec,t)-f(\bfx,\hat{\thetavec},t)-\varphi(\psi,\omegavec,t)+\varepsilon(t),
\end{equation}
where, if not stated overwise, the function
$\varepsilon:\Real_+\rightarrow\Real$, $\varepsilon\in L_{2}^1
[t_0,\infty]\cap C^0$. One of the immediate advantages of equation
(\ref{error_model_d}) in comparison with (\ref{error_model}) is
that it allows us to take the presence of coupling between
interconnected systems into consideration.

Let us now specify the desired properties of the function
$\varphi(\psi,\omegavec,t)$ in (\ref{control}),
(\ref{error_model_d}). The majority of known algorithms for
parameter estimation and adaptive control
\cite{Kokotovich95,Fradkov99,Narendra89,Sastry89} assume global
(Lyapunov) stability
 of system
(\ref{error_model_d}) for $\thetavec\equiv\hat{\thetavec}$.  In
our study, however, we refrain from this standard, restrictive
requirement. Instead we propose that finite energy of the signal
$f(\bfx(t),\thetavec,t)-f(\bfx(t),\hat{\thetavec}(t),t)$, defined
for example by its $L_{2}^1[t_0,\infty]$ norm with respect to the
variable $t$, results in finite deviation from the target set
given by the equality $\psi(\bfx,t)=0$. Formally this requirement
is introduced in Assumption \ref{assume:gain}:
\begin{assume}[Target dynamics operator]\label{assume:gain} Consider the following system:
\begin{equation}\label{eq:target_dynamics}
\dpsi=-\varphi(\psi,\omegavec,t)+\zeta(t),
\end{equation}
where $\zeta:\Real_+\rightarrow\Real$ and
$\varphi(\psi,\omegavec,t)$ is defined in (\ref{error_model_d}).
Then for every $\omegavec\in\Omega_\omega$ system
(\ref{eq:target_dynamics}) has $L_{2}^1 [t_0,\infty]\mapsto
L_\infty^1[t_0,\infty]$ gain with respect to input $\zeta(t)$. In
other words, there exists a function $\gamma_{\infty,2}$ such that
\begin{equation}\label{eq:gain_psi_L2}
\|\psi(t)\|_{\infty,[t_0,T]}\leq
\gamma_{\infty,2}(\psi_0,\omegavec,\|\zeta(t)\|_{2,[t_0,T]}), \ \
\forall \ \zeta(t)\in L_{2}^1[t_0,T]
\end{equation}
\end{assume}
In contrast to conventional approaches, Assumption
\ref{assume:gain} does not require global {\it asymptotic
stability} of the origin of the unperturbed (i.e for $\zeta(t)=0$)
system (\ref{eq:target_dynamics}). When the stability of the
target dynamics $\dpsi=-\varphi(\psi,\omegavec,t)$ is known
a-priori, one of the benefits of Assumption \ref{assume:gain} is
that there is no need to know a {\it particular Lyapunov function}
of the unperturbed system.

So far we have introduced basic assumptions on system
(\ref{system1})  and the class of feedback considered in this
article. Let us now specify the class of functions
$f(\bfx,\thetavec,t)$ in (\ref{error_model_d}). Since general
parametrization of function $f(\bfx,\thetavec,t)$ is
methodologically difficult to deal with, but solutions provided
for nonlinearities with convenient linear re-parametrization often
yield physically implausible models and large number of unknown
parameters, we have opted for a new class of parameterizations.
As a candidate for such a
parametrization we suggest nonlinear functions that satisfy the
following assumption:
\begin{assume}[Monotonicity and Growth Rate in
Parameters]\label{assume:alpha}For the given function
$f(\bfx,\thetavec,t)$ in (\ref{error_model_d}) there exists
function  $\alphavec(\bfx,t): \Real^{n}\times \Real_+\rightarrow
\Real^d, \ \alphavec(\bfx,t)\in \mathcal{C}^1$ and positive
constant $D>0$ such that
\begin{equation}\label{eq:assume_alpha}
(f(\bfx,\hat{\thetavec},t)-f(\bfx,\thetavec,t))(\alphavec(\bfx,t)^{T}(\hat{\thetavec}-\thetavec))\geq0
\end{equation}
\begin{equation}\label{eq:assume_gamma}
|f(\bfx,\hat{\thetavec},t)-f(\bfx,\thetavec,t)|\leq D
|\alphavec(\bfx,t)^{T}(\hat{\thetavec}-\thetavec)|
\end{equation}
\end{assume}
This set of conditions naturally extends  from systems that are
linear in parameters  to those with nonlinear parametrization.
Examples and models of physical and artificial systems which
satisfy Assumption \ref{assume:alpha} (at least for bounded
$\thetavec,\hat{\thetavec}\in \Omega_\theta$) can be found in the
following references
\cite{Armstrong_1993,Boskovic_1995,Canudas_1999,Abbott_2001,Kitching_2000}.
Assumption \ref{assume:alpha} bounds the growth rate of the
difference $|f(\bfx,\thetavec,t)-f(\bfx,\hat{\thetavec},t)|$ by
the functional
$D|\alphavec(\bfx,t)^{T}(\hat{\thetavec}-\thetavec)|$.
In addition, it might also be useful to have an estimate of
$|f(\bfx,\thetavec,t)-f(\bfx,\hat{\thetavec},t)|$ from below, as
specified in Assumption \ref{assume:alpha_upper}:
\begin{assume}\label{assume:alpha_upper} For the given function
$f(\bfx,\thetavec,t)$ in (\ref{error_model_d}) and function
$\alphavec(\bfx,t)$, satisfying Assumption \ref{assume:alpha},
there exists a positive constant $D_1>0$  such that
\begin{equation}\label{eq:assume_alpha_upper}
|f(\bfx,\hat{\thetavec},t)-f(\bfx,\thetavec,t)|\geq D_1
|\alphavec(\bfx,t)^{T}(\hat{\thetavec}-\thetavec)|
\end{equation}
\end{assume}

\noindent In problems of adaptation, parameter and optimization
estimation, effectiveness of the algorithms often depends on how
"good" the nonlinearity $f(\bfx,\thetavec,t)$ is, and how
predictable is the system's behavior. As a measure of goodness and
predictability usually the substitutes as smoothness and
boundedness are considered. In our study, we distinguish several
of such specific properties of the functions $f(\bfx,\thetavec,t)$
and $\varphi(\psi,\omegavec,t)$. These properties are provided
below.

\begin{hyp}\label{hyp:locally_bound_uniform_f} The function $f(\bfx,\thetavec,t)$ is locally bounded
with respect to $\bfx$, ${\thetavec}$  uniformly in $t$.
\end{hyp}

\begin{hyp}\label{hyp:locally_bound_uniform_df} The function $f(\bfx,\thetavec,t)\in \mathcal{C}^1$, and $ \pd
{f(\bfx,\thetavec,t)}/{\pd t}$ is locally bounded with respect to
$\bfx$, ${\thetavec}$ uniformly in $t$.
\end{hyp}

\begin{hyp}\label{hyp:locally_bound_uniform_phi} The function $\varphi(\psi,\omegavec,t)$ is locally bounded in $\psi$,
$\omegavec$ uniformly in $t$.
\end{hyp}

Let us show that under an additional structural requirement, which
relates properties of the function $\alphavec(\bfx,t)$ and
vector-field
$\bff(\bfx,\thetavec)=\bff_1(\bfx,\thetavec)\oplus\bff_2(\bfx,\thetavec)$
in (\ref{system1}), (\ref{system}), there exist adaptive
algorithms ensuring that the following desired property holds:
\begin{equation}\label{eq:desired_prop}
\bfx(t)\in L_\infty^n[t_0,\infty]; \
f(\bfx(t),\thetavec,t)-f(\bfx,\hat{\thetavec}(t),t)\in
L_{2}^1[t_0,\infty]
\end{equation}

Consider the following adaptation algorithms:
\begin{equation}\label{fin_forms_ours_tr1}
\begin{split}
\hat{\thetavec}(\bfx,t)&=\Gamma(\hat{\thetavec}_P(\bfx,t)+\hat{\thetavec}_I(t));
\  \Gamma\in\Real^{d\times d}, \ \Gamma>0
\\ \hat{\thetavec}_P(\bfx,t)&=
\psi(\bfx,t)\alphavec(\bfx,t)-\Psi(\bfx,t) \\
\dot{\hat{\thetavec}}_I&=\varphi(\psi(\bfx,t),\omegavec,t)\alphavec(\bfx,t)+\mathcal{R}(\bfx,\hat{\thetavec},u(\bfx,\hat{\thetavec},t),t),
\end{split}
\end{equation}
where the function
$\mathcal{R}(\bfx,\hat{\thetavec},u(\bfx,\hat{\thetavec},t),t):\Real^n\times\Real^d\times\Real\times\Real_+\rightarrow\Real^d$
in (\ref{fin_forms_ours_tr1}) is given as follows:
\begin{equation}\label{fin_forms_ours_tr11}
\begin{split}
&\mathcal{R}(\bfx,u(\bfx,\hat{\thetavec},t),t)={\pd
\Psi(\bfx,t)}/{\pd t}-\psi(\bfx,t)({\pd
\alphavec(\bfx,t)}/{\pd t}+L_{\bff_1}\alphavec(\bfx,t))\\
&    + L_{\bff_1}
\Psi(\bfx,t)-(\psi(\bfx,t)L_{\bfg_1}\alphavec(\bfx,t)-L_{\bfg_1}
\Psi(\bfx,t))u(\bfx,\hat{\thetavec},t)
\end{split}
\end{equation}
and function
$\Psi(\bfx,t):\Real^{n}\times\Real_+\rightarrow\Real_d$,
$\Psi(\bfx,t)\in \mathcal{C}^1$  satisfies Assumption
\ref{assume:explicit_realizability}.
\begin{assume}\label{assume:explicit_realizability} There exists a function $\Psi(\bfx,t)$ such that
\begin{equation}\label{eq:assume_explicit}
\frac{\pd \Psi(\bfx,t)}{\pd \bfx_2}-\psi(\bfx,t)\frac{\pd
\alphavec(\bfx,t)}{\pd \bfx_2}=0
\end{equation}
\end{assume}
Additional restrictions imposed by this assumption will be
discussed in some details after we summarize the properties of
system (\ref{system1}), (\ref{control}),
(\ref{fin_forms_ours_tr1}), (\ref{fin_forms_ours_tr11}) in the
following theorem.

\begin{theorem}[Properties of the decoupled systems]\label{stability_theorem}
Let system (\ref{system1}), (\ref{error_model_d}),
(\ref{fin_forms_ours_tr1}), (\ref{fin_forms_ours_tr11}) be given
and Assumptions \ref{assume:alpha}, \ref{assume:alpha_upper},
\ref{assume:explicit_realizability} be satisfied. Then the
following properties hold

P1) Let for the given initial conditions $\bfx(t_0)$,
$\hat{\thetavec}_I(t_0)$ and parameters vector $\thetavec$,
interval $[t_0,T^\ast]$ be the (maximal) time-interval of
existence of solutions of the closed loop system (\ref{system1}),
(\ref{error_model_d}), (\ref{fin_forms_ours_tr1}),
(\ref{fin_forms_ours_tr11}). Then
\begin{equation}\label{eq:f_diff_L2}
\|f(\bfx(t),\thetavec,t)-f(\bfx(t),\hat{\thetavec}(t),t))\|_{2,[t_0,T^\ast]}\leq
D_f(\thetavec,t_0,\Gamma,\|\varepsilon(t)\|_{2,[t_0,T^\ast]});
\end{equation}
\[
D_f(\thetavec,t_0,\Gamma,\|\varepsilon(t)\|_{2,[t_0,T^\ast]})=\left(\frac{D}{2}\|\thetavec-\hat{\thetavec}(t_0)\|^{2}_{\Gamma^{-1}}\right)^{0.5}
+ \frac{D}{D_1}\|\varepsilon(t)\|_{2,[t_0,T^\ast]}
\]
\[
\|\thetavec-\hat\thetavec(t)\|^{2}_{\Gamma^{-1}}\leq
\|\hat{\thetavec}(t_0)-\thetavec\|^{2}_{\Gamma^{-1}}+\frac{D}{2
D_1^2}\|\varepsilon(t)\|^{2}_{2,[t_0,T^\ast]}
\]

\noindent  In addition, if Assumptions \ref{assume:psi} and
\ref{assume:gain} are satisfied then

P2)  $\psi(\bfx(t),t)\in L_\infty^1[t_0,\infty]$, $\bfx(t)\in
L_{\infty}^n[t_0,\infty]$ and
\begin{equation}\label{eq:psi_gain}
\|\psi(\bfx(t),t)\|_{\infty,[t_0,\infty]}\leq
\gamma_{\infty,2}\left(\psi(\bfx_0,t_0),\omegavec,\mathcal{D}\right)
\end{equation}
\[
\mathcal{D}=D_f(\thetavec,t_0,\Gamma,\|\varepsilon(t)\|_{2,[t_0,\infty]})+\|\varepsilon(t)\|_{2,[t_0,\infty]}
\]

P3) if properties  H\ref{hyp:locally_bound_uniform_f},
H\ref{hyp:locally_bound_uniform_phi} hold, and system
(\ref{eq:target_dynamics}) has $L_{2}^1 [t_0,\infty]\mapsto
L_{p}^1 [t_0,\infty]$, $p>1$ gain with respect to input $\zeta(t)$
and output $\psi$ then
\begin{equation}\label{eq:convergence_psi_theorem}
\varepsilon(t)\in L_{2}^1 [t_0,\infty]\cap
L_{\infty}^1[t_0,\infty]\Rightarrow
\lim_{t\rightarrow\infty}\psi(\bfx(t),t)=0
\end{equation}

If, in addition, property H\ref{hyp:locally_bound_uniform_df}
holds, and the functions $\alphavec(\bfx,t)$, $\pd
\psi(\bfx,t)/\pd t$ are locally bounded with respect to $\bfx$
uniformly in $t$, then

P4) the following holds
\begin{equation}\label{eq:convergence_f_theorem}
\lim_{t\rightarrow\infty}f(\bfx(t),\thetavec,t)-f(\bfx(t),\hat{\thetavec}(t),t)=0
\end{equation}

\end{theorem}
The proof of Theorem \ref{stability_theorem} and subsequent
results are given in Section 6.

Let us briefly comment on Assumption
\ref{assume:explicit_realizability}.
 Let $\alphavec(\bfx,t)\in
\mathcal{C}^2$,
$\alphavec(\bfx,t)=\mathrm{col}(\alpha_1(\bfx,t),\dots,\alpha_d(\bfx,t))$,
then necessary  and sufficient conditions for existence of the
function $\Psi(\bfx,t)$ follow from the
Poincar$\acute{\mathrm{e}}$ lemma:
\begin{equation}\label{eq:poincare}
\frac{\pd}{\pd \bfx_2}\left(\psi(\bfx,t)\frac{\pd
\alpha_i(\bfx,t)}{\pd
\bfx_2}
\right)=\left(\frac{\pd}{\pd
\bfx_2}\left(\psi(\bfx,t)\frac{\pd \alpha_i(\bfx,t)}{\pd
\bfx_2}
\right)
\right)^T
\end{equation}
This relation, in the form of conditions of existence of the
solutions for function $\Psi(\bfx,t)$ in
(\ref{eq:assume_explicit}), takes into account structural
properties of system (\ref{system1}), (\ref{error_model_d}).
Indeed,
consider partial derivatives $\pd \alpha_i(\bfx,t)/\pd \bfx_2$,
$\pd \psi(\bfx,t)/\pd \bfx_2$ with respect to the vector
$\bfx_2=(x_{21},\dots,x_{2p})^T$. Let
\begin{equation}\label{eq:single_dim}
\begin{split}
\frac{\pd \psi(\bfx,t)}{\pd \bfx_2}=\left(\begin{array}{cccccccc}
                                     0& 0
                                     & \cdots & 0& \ast & 0&\cdots&0
                                    \end{array}\right), \
\frac{\pd
\alpha_i(\bfx,t)}{\pd\bfx_2}=\left(\begin{array}{cccccccc}
                                     0 & 0
                                     & \cdots & 0&
                                     \ast &
                                     0&\cdots&0
                                        \end{array}\right)
\end{split}
\end{equation}
where the symbol $\ast$ denotes a function of $\bfx$ and $t$. Then
condition (\ref{eq:single_dim}) guarantees that equality
(\ref{eq:poincare}) (and, subsequently, Assumption
\ref{assume:explicit_realizability}) holds. In case $\pd
\alpha(\bfx_1\oplus \bfx_2,t)/\pd \bfx_2=0$, Assumption
\ref{assume:explicit_realizability} holds for arbitrary
$\psi(\bfx,t)\in \mathcal{C}^1$. If $\psi(\bfx,t)$,
$\alphavec(\bfx,t)$ depend on a single component of $\bfx_2$, for
instance $x_{2k}, \ k\in\{0,\dots,p\}$, then conditions
(\ref{eq:single_dim}) hold and the function $\Psi(\bfx,t)$ can be
derived explicitly by integration
\begin{equation}\label{eq:single_dim_int}
\Psi(\bfx,t)=\int\psi(\bfx,t)\frac{\alphavec(\bfx,t)}{\pd x_{2k}}d
x_{2k}
\end{equation}
In all other cases, existence of the required function
$\Psi(\bfx,t)$ follows from (\ref{eq:poincare}).

In the general case, when $\dim\{\bfx_2\}>1$, the problems of
finding a function $\Psi(\bfx,t)$ satisfying condition
(\ref{eq:assume_explicit}) can be avoided (or converted into one
with an already known solutions such as (\ref{eq:poincare}),
(\ref{eq:single_dim_int})) by the {\it embedding} technique
proposed in \cite{ECC_2003}. The main idea of the method is to
introduce an auxiliary system that is forward-complete with
respect to input $\bfx(t)$
\begin{equation}\label{eq:embed}
\begin{split}
\dot{\xivec}&=\bff_\xivec(\bfx,\xivec,t), \ \xivec\in\Real^z \\
\bfh_\xi&=\bfh_\xi(\xivec,t), \
\Real^z\times\Real_+\rightarrow\Real^h
\end{split}
\end{equation}
such that
\begin{equation}\label{eq:embed_L2}
\|f(\bfx(t),\thetavec,t)-f(\bfx_1(t)\oplus\bfh_\xi(t)\oplus\bfx_2'(t),\thetavec,t)\|_{2,[t_0,T]}
\leq C_\xi\in\Real_+
\end{equation}
for all $T\geq t_0$, and $\dim\{{\bfh_\xi}\}+\dim{\{\bfx_2'\}}=p$.
Then (\ref{error_model_d}) can be rewritten as follows:
\begin{equation}\label{error_model_d1}
\dpsi=f(\bfx_1\oplus\bfh_\xi\oplus\bfx_2',\thetavec,t)-f(\bfx_1\oplus\bfh_\xi\oplus\bfx_2',\hat\thetavec,t)-\varphi(\psi,\omegavec,t)+\varepsilon_\xi(t),
\end{equation}
where $\varepsilon_\xi(t)\in L_{2}^1 [t_0,\infty]$, and
$\dim\{\bfx_2'\}=p-h<p$. In principle, the dimension of $\bfx_2'$
could be reduced to $1$ or $0$. As soon as this is ensured,
Assumption \ref{assume:explicit_realizability} will be satisfied
and the results of Theorem \ref{stability_theorem} follow.
Sufficient conditions ensuring the existence of such an embedding
in the general case are provided in \cite{ECC_2003}. For systems
in which the parametric uncertainty  can be reduced to vector
fields with low-triangular structure the embedding is given in
\cite{ALCOSP_2004}.


\section{Main Results}\label{sec:main}

Without loss of generality  let us rewrite interconnection
(\ref{eq:system:s1}), (\ref{eq:system:s2}) as
follows
:
\begin{equation}\label{eq:system:s11}
\begin{split}
\dot{\bfx}_1&=\bff_1(\bfx)+\bfg_1(\bfx)u_x\\
\dot{\bfx}_2 &=\bff_2(\bfx,\thetavec_x)+\gamma_y(\bfy,t)+
\bfg_2(\bfx)u_x
\end{split}
\end{equation}

\begin{equation}\label{eq:system:s21}
\begin{split}
\dot{\bfy}_1&=\bfq_1(\bfy)+\bfz_1(\bfy)u_y\\
\dot{\bfy}_2&=\bfq_2(\bfy,\thetavec_y)+\gamma_x(\bfx,t)+\bfz_2(\bfy)u_y
\end{split}
\end{equation}

Let us now consider the following control functions
\begin{equation}\label{control_s1}
\begin{split}
u_x(\bfx,\hat{\thetavec}_x,\omegavec_x,t)&=(L_{\bfg(\bfx)}\psi_x(\bfx,t))^{-1}\left(-L_{\bff(\bfx,\hat{\thetavec}_x)}\psi_x(\bfx,t)-\varphi_x(\psi_x,\omegavec_x,t)\right.\\
& \left.-\frac{\pd\psi_x(\bfx,t)}{\pd t}\right), \ \ \varphi_x: \
\Real\times\Real^w\times\Real_+\rightarrow\Real
\end{split}
\end{equation}
\begin{equation}\label{control_s2}
\begin{split}
u_y(\bfy,\hat{\thetavec}_y,\omegavec_y,t)&=(L_{\bfz(\bfy)}\psi_y(\bfy,t))^{-1}\left(-L_{\bfq(\bfy,\hat{\thetavec}_y)}\psi_y(\bfy,t)-\varphi_y(\psi_y,\omegavec_y,t)\right.\\
&\left.-\frac{\pd\psi_y(\bfy,t)}{\pd t}\right), \ \ \varphi_y: \
\Real\times\Real^w\times\Real_+\rightarrow\Real
\end{split}
\end{equation}
These functions transform the original equations
(\ref{eq:system:s11}), (\ref{eq:system:s21}) into the following
form
\begin{equation}\label{eq:error_coupled}
\begin{split}
\dpsi_x&=-\varphi_x(\psi_x,\omegavec_x,t)+f_x(\bfx,\thetavec_x,t)-f_x(\bfx,\hat{\thetavec}_x,t)+h_y(\bfx,\bfy,t)\\
\dpsi_y&=-\varphi_y(\psi_x,\omegavec_y,t)+f_y(\bfy,\thetavec_y,t)-f_y(\bfy,\hat{\thetavec}_y,t)+h_x(\bfx,\bfy,t),
\end{split}
\end{equation}
where
\[
h_x(\bfx,\bfy,t)=L_{\gamma_y(\bfy,t)}\psi_x(\bfx,t), \
h_y(\bfx,\bfy,t)=L_{\gamma_x(\bfx,t)}\psi_y(\bfy,t)
\]
\[
f_x(\bfx,\thetavec_x,t)=L_{\bff(\bfx,\thetavec_x)}\psi_x(\bfx,t), \
f_y(\bfx,\thetavec_y,t)=L_{\bfq(\bfy,\thetavec_y)}\psi_y(\bfy,t)
\]

Consider the following adaptation algorithms
\begin{equation}\label{fin_forms_ours_tr1x}
\begin{split}
\hat{\thetavec}_x(\bfx,t)&=\Gamma_x(\hat{\thetavec}_{P,x}(\bfx,t)+\hat{\thetavec}_{I,x}(t));
\  \Gamma_x\in\Real^{d\times d}, \ \Gamma_x>0
\\ \hat{\thetavec}_{P,x}(\bfx,t)&=
\psi_x(\bfx,t)\alphavec_x(\bfx,t)-\Psi_x(\bfx,t) \\
\dot{\hat{\thetavec}}_{I,x}&=\varphi_x(\psi_x(\bfx,t),\omegavec_x,t)\alphavec_x(\bfx,t)+\mathcal{R}_x(\bfx,\hat{\thetavec}_x,u_x(\bfx,\hat{\thetavec}_x,t),t),
\end{split}
\end{equation}

\begin{equation}\label{fin_forms_ours_tr1y}
\begin{split}
\hat{\thetavec}_y(\bfx,t)&=\Gamma_y(\hat{\thetavec}_{P,y}(\bfy,t)+\hat{\thetavec}_{I,y}(t));
\  \Gamma_y\in\Real^{d\times d}, \ \Gamma_y>0
\\ \hat{\thetavec}_{P,y}(\bfy,t)&=
\psi_y(\bfy,t)\alphavec_y(\bfy,t)-\Psi_y(\bfy,t) \\
\dot{\hat{\thetavec}}_{I,y}&=\varphi_y(\psi_y(\bfy,t),\omegavec_y,t)\alphavec_y(\bfy,t)+\mathcal{R}_y(\bfx,\hat{\thetavec}_y,u_y(\bfy,\hat{\thetavec}_y,t),t),
\end{split}
\end{equation}
where $\mathcal{R}_x(\cdot)$, $\mathcal{R}_y(\cdot)$ are defined
as in (\ref{fin_forms_ours_tr11}), and the functions
$\Psi_x(\cdot)$, $\Psi_y(\cdot)$ will be specified later. Now we
are ready to formulate the following result

\begin{theorem}[Properties of the interconnected systems]\label{theorem:interconnection} Let systems (\ref{eq:system:s11}), (\ref{eq:system:s21}) be
given. Furthermore, suppose that the following conditions hold:

1) The functions $\psi_x(\bfx,t)$, $\psi_y(\bfy,t)$ satisfy
Assumption \ref{assume:psi} for systems (\ref{eq:system:s11}),
(\ref{eq:system:s21}) respectively;

2) The systems
\begin{equation}\label{eq:target_dynamics_connected}
\dot{\psi}_x=-\varphi_x(\psi_x,\omegavec_x,t)+\zeta_x(t), \ \
\dot{\psi}_y=-\varphi_y(\psi_y,\omegavec_y,t)+\zeta_y(t)
\end{equation}
satisfy Assumption \ref{assume:gain} with corresponding mappings
\[
\gamma_{x_{\infty,2}}(\psi_{x_0},\omegavec_x,\|\zeta_x(t)\|_{2,[t_0,T]}),
\ \
\gamma_{y_{\infty,2}}(\psi_{y_0},\omegavec_y,\|\zeta_y(t)\|_{2,[t_0,T]}),
\]

3) The systems (\ref{eq:target_dynamics_connected}) have
$L_2^1[t_0,\infty]\mapsto L_2^1[t_0,\infty]$ gains, that is
\begin{equation}\label{eq:L_2_2_gains}
\begin{split}
\|\psi_x(\bfx(t),t)\|_{2,[t_0,T]}&\leq
C_{\gamma_x}+\gamma_{x_{2,2}}(\|\zeta_x(t)\|_{2,[t_0,T]}),\\
\|\psi_y(\bfy(t),t)\|_{2,[t_0,T]}&\leq
C_{\gamma_y}+\gamma_{y_{2,2}}(\|\zeta_y(t)\|_{2,[t_0,T]}),\\
C_{\gamma_x}, \ C_{\gamma_y}\in\Real_+& \gamma_{x_{2,2}}, \
\gamma_{y_{2,2}}\in\mathcal{K}_\infty
\end{split}
\end{equation}

4) The functions $f_x(\bfx,\thetavec_x,t)$,
$f_y(\bfy,\thetavec_y,t)$ satisfy Assumptions \ref{assume:alpha},
\ref{assume:alpha_upper} with corresponding constants $D_x$,
$D_{x_1}$, $D_y$, $D_{y_1}$ and functions $\alphavec_x(\bfx,t)$,
$\alphavec_y(\bfy,t)$;

5) The functions $h_x(\bfx,\bfy,t)$, $h_y(\bfx,\bfy,t)$ satisfy
the following inequalities:
\begin{equation}\label{eq:disturbance_gain}
\|h_x(\bfx,\bfy,t)\|\leq \beta_x \|\psi_x(\bfx,t)\|, \
\|h_y(\bfx,\bfy,t)\|\leq \beta_y \|\psi_y(\bfy,t)\|, \  \beta_x,
\beta_y\in \Real_+
\end{equation}

Finally, let the functions $\Psi_x(\bfx,t)$, $\Psi_y(\bfy,t)$ in
(\ref{fin_forms_ours_tr1x}), (\ref{fin_forms_ours_tr1y}) satisfy
Assumption \ref{assume:explicit_realizability}
for systems (\ref{eq:system:s11}), (\ref{eq:system:s21})
respectively, and there exist functions $\rho_1(\cdot), \
\rho_2(\cdot), \ \rho_3(\cdot)>Id(\cdot)\in\mathcal{K}_\infty$ and
constant $\bar{\Delta}\in\Real_+$ such the following inequality
holds:
\begin{equation}\label{eq:small_gain_adapt}
\beta_y\circ\gamma_{y_{2,2}}\circ\rho_1\circ\left(\frac{D_y}{D_{y,1}}+1\right)\circ\rho_3\circ
\beta_x\circ
\gamma_{x_{2,2}}\circ\rho_2\circ\left(\frac{D_x}{D_{x,1}}+1\right)(\Delta)<
\Delta
\end{equation}
for all $\Delta\geq \bar{\Delta}$. Then

C1) The interconnection (\ref{eq:system:s11}),
(\ref{eq:system:s21}) with controls (\ref{control_s1}),
(\ref{control_s2}) is forward-complete and trajectories $\bfx(t)$,
$\bfy(t)$ are bounded

Furthermore,

C2) if properties  H\ref{hyp:locally_bound_uniform_f},
H\ref{hyp:locally_bound_uniform_phi} hold for
$f_x(\bfx,\thetavec_x,t)$,  $f_y(\bfy,\thetavec_y,t)$,
$h_x(\bfx,\bfy,t)$, $h_y(\bfx,\bfy,t)$,  and also functions
$\varphi_x(\psi_x,\omegavec_x,t)$,
$\varphi_y(\psi_y,\omegavec_y,t)$, then
\begin{equation}\label{eq:convergence_psi_xy}
\lim_{t\rightarrow\infty}\psi_x(\bfx(t),t)=0, \
\lim_{t\rightarrow\infty}\psi_y(\bfy(t),t)=0
\end{equation}

Moreover,

C3) if property H\ref{hyp:locally_bound_uniform_df} holds for
$f_x(\bfx,\thetavec_x,t)$, $f_y(\bfy,\thetavec_y,t)$, and the
functions
\[
\alphavec_x(\bfx,t), \ \pd \psi_x(\bfx,t)/\pd t, \
\alphavec_y(\bfy,t), \ \pd \psi_y(\bfy,t)/\pd t
\]
are locally bounded with respect to $\bfx$, $\bfy$ uniformly in $t$,
then
\begin{equation}\label{eq:convergence_f_xy}
\begin{split}
\lim_{t\rightarrow\infty}f_x(\bfx(t),\thetavec_x,t)-f_x(\bfx(t),\hat{\thetavec}_x(t),t)&=0,
\\
\lim_{t\rightarrow\infty}f_y(\bfy(t),\thetavec_y,t)-f_y(\bfy(t),\hat{\thetavec}_y(t),t)&=0
\end{split}
\end{equation}
\end{theorem}

Let us briefly comment on the conditions and assumptions of
Theorem \ref{theorem:interconnection}. Conditions 1), 2) specify
restrictions on the goal functionals, similar to those of Theorem
\ref{stability_theorem}. Condition 3) is analogous to requirement
to P3) in Theorem \ref{stability_theorem}, condition 5) specifies
uncertainties in the coupling functions $h_x(\cdot)$, $h_y(\cdot)$
in terms of their growth rates w.r.t. $\psi_x(\cdot)$,
$\psi_y(\cdot)$. We observe here that this property is needed in
order to characterize the $L_2$ norms of functions
$h_x(\bfx(t),\bfy(t),t)$, $h_y(\bfx(t),\bfy(t),t)$ in terms of the
$L_2$ norms of functions $\psi_x(\bfx(t),t)$, $\psi_y(\bfy(t),t)$.
Therefore, it is possible to replace requirement
(\ref{eq:disturbance_gain}) with the following set of conditions:
\begin{equation}\label{eq:disturbance_gain_1}
\begin{split}
 \|h_x(\bfx(t),\bfy(t),t)\|_{2,[t_0,T]}&\leq \beta_x
\|\psi_x(\bfx(t),t)\|_{2,[t_0,T]}+C_x, \\
 \|h_y(\bfx(t),\bfy(t),t)\|_{2,[t_0,T]}&\leq \beta_y
\|\psi_y(\bfy(t),t)\|_{2,[t_0,T]}+C_y
\end{split}
\end{equation}
The replacement will allow us to extend results of Theorem
\ref{theorem:interconnection} to interconnections of systems where
the coupling functions do not depend explicitly on
$\psi_x(\bfx(t),t)$, $\psi_y(\bfy(t),t)$. We illustrate this
possibility later with an example.

Condition (\ref{eq:small_gain_adapt}) is the small-gain condition
with respect to the $L_2^1[t_0,T]$ norms for interconnection
(\ref{eq:system:s11}), (\ref{eq:system:s21}) with control
(\ref{control_s1}), (\ref{control_s2}). In the case that mappings
$\gamma_{x_{2,2}}(\cdot)$, $\gamma_{y_{2,2}}(\cdot)$ in
(\ref{eq:target_dynamics_connected}) are majorated by linear
functions
\[
\gamma_{x_{2,2}}(\Delta)\leq g_{x_{2,2}} \Delta, \
\gamma_{y_{2,2}}(\Delta)\leq g_{y_{2,2}} \Delta, \ \Delta\geq 0,
\]
condition (\ref{eq:small_gain_adapt}) reduces to the much simpler
\[
\beta_y \beta_x g_{x_{2,2}} g_{y_{2,2}}
\left(\frac{D_y}{D_{y,1}}+1\right)\left(\frac{D_x}{D_{x,1}}+1\right)<
1
\]
Notice also that the mappings $\gamma_{x_{2,2}}(\cdot)$,
$\gamma_{y_{2,2}}(\cdot)$ are defined by properties of the target
dynamics (\ref{eq:target_dynamics_connected}), and, in principle,
these can be made arbitrarily small. This eventually leads to the
following conclusion: the smaller the $L_2$-gains of the target
dynamics of systems $\mathcal{S}_1$, $\mathcal{S}_2$, the wider
the class of
 nonlinearities (bounds for $\beta_x$, $\beta_y$, domains
of $D_x$, $D_{1,x}$, $D_y$, $D_{1,y}$) which admit a solution to
Problem \ref{problem:decentralized}.

\paragraph{Example}

Let us illustrate application of Theorem
\ref{theorem:interconnection} to the problem of decentralized
control of two coupled oscillators with nonlinear damping.
Consider the following interconnected systems:
\begin{equation}\label{eq:example_dec_model}
\left\{\begin{array}{ll}
\dot{x}_{1}&=x_{2}\\
\dot{x}_{2}&=f_x(x_{1},\theta_x)+k_1 y_{1} + u_x,
\end{array} \right. \ \
\left\{
\begin{array}{ll}
\dot{y}_{1}&=y_{2}\\
\dot{y}_{22}&=f_y(y_{1},\theta_y)+k_2 x_{1}+ u_y,
\end{array}\right.
\end{equation}
where $k_1$, $k_2\in\Real$ are uncertain parameters of coupling,
functions $f(x_{1},\theta_x)$, $f(y_{1},\theta_y)$
 stand for the nonlinear damping terms, and
$\theta_{x}$, $\theta_y$ are unknown parameters. For illustrative
purpose we assume the following mathematical model for functions
$f_x(\cdot)$, $f_y(\cdot)$ in (\ref{eq:example_dec_model}):
\begin{equation}\label{eq:example_dec_uncertainty}
\begin{split}
f_x(x_{1},\theta_x)&= \theta_x (x_{1}-x_0)+0.5\sin
(\theta_x(x_{1}-x_0)),\\
\ f_y(y_{1},\theta_y)&= \theta_y (y_{1}-y_0)+0.6\sin
(\theta_y(y_{1}-y_0))
\end{split}
\end{equation}
where $x_0$, $y_0$ are known. Let the control goal be to steer
states $\bfx$ and $\bfy$ to the origin. Consider the following
goal functions
\begin{equation}\label{eq:example_psi}
\psi_x(\bfx,t)=x_1+x_2, \ \psi_y(\bfy,t)= y_1+y_2
\end{equation}
Taking into account equations (\ref{eq:example_dec_model}) and
(\ref{eq:example_psi}) we can derive that
\begin{equation}\label{eq:example_relative_dynamics}
\dot{x}_1=-x_1+\psi_x(\bfx(t),t), \ \dot{y}_1=-y_1+\psi_y(\bfy,t)
\end{equation}
This automatically implies that
\[
\begin{split}
\|x_1(t)\|_{\infty,[t_0,T]}&\leq
\|x_1(t_0)\|+\|\psi_x(\bfx(t),t)\|_{\infty,[t_0,T]}\\
\|y_1(t)\|_{\infty,[t_0,T]}&\leq
\|y_1(t_0)\|+\|\psi_y(\bfy(t),t)\|_{\infty,[t_0,T]}
\end{split}
\]
Hence, Assumption \ref{assume:psi} is satisfied for chosen goal
functions $\psi_x(\cdot)$ and $\psi_y(\cdot)$.  Notice also that
equalities (\ref{eq:example_relative_dynamics}) imply that
\begin{equation}\label{eq:example_L2_gains}
\begin{split}
\|x_1(t)\|_{2,[t_0,T]}&\leq 2^{-1/2}\|x_1(t_0)\|+
\|\psi_x(\bfx,t)\|_{2,[t_0,T]}\\
\|y_1(t)\|_{2,[t_0,T]}&\leq 2^{-1/2}\|y_1(t_0)\|+
\|\psi_y(\bfy,t)\|_{2,[t_0,T]}
\end{split}
\end{equation}
Moreover, according to (\ref{eq:example_relative_dynamics})
limiting relations
\begin{equation}\label{eq:example_control_goal_limit}
\begin{split}
&
\lim_{t\rightarrow\infty}\psi_x(\bfx(t),t)=\lim_{t\rightarrow\infty}x_1(t)+x_2(t)=0,\\
&
\lim_{t\rightarrow\infty}\psi_y(\bfy(t),t)=\lim_{t\rightarrow\infty}y_1(t)+y_2(t)=0
\end{split}
\end{equation}
guarantee that
\[
\lim_{t\rightarrow\infty} x_1(t)=0, \
\lim_{t\rightarrow\infty}x_2(t)=0, \ \lim_{t\rightarrow\infty}
y_1(t)=0, \ \lim_{t\rightarrow\infty}y_2(t)=0
\]
Hence, property (\ref{eq:example_control_goal_limit}) ensures
asymptotic reaching of the control goal.

According to equations (\ref{control_s1}), (\ref{control_s2})
control functions
\begin{equation}\label{eq:example_control}
\begin{split}
u_x&=-\lambda_x\psi_x-x_2-f_x(x_1,\hat{\theta}_x)\\
u_y&=-\lambda_y\psi_y-y_2-f_y(y_1,\hat{\theta}_y), \ \lambda_x, \
\lambda_y>0
\end{split}
\end{equation}
transform system (\ref{eq:example_dec_model}) into the following
form
\begin{equation}\label{eq:example_error_model}
\begin{split}
\dot{\psi}_x&=-\lambda_x \psi_x +
f_x(x_1,\theta_x)-f_x(x_1,\hat{\theta}_x)+k_1 y_1\\
\dot{\psi}_x&=-\lambda_x \psi_x +
f_x(x_1,\theta_x)-f_x(x_1,\hat{\theta}_x)+k_2 x_1
\end{split}
\end{equation}
Notice that systems
\[
\dot{\psi}_x=-\lambda_x \psi_x +\xi_x(t), \
\dot{\psi}_y=-\lambda_y \psi_t +\xi_y(t)
\]
satisfy Assumption \ref{assume:gain} with
\[
\gamma_{x_{2,2}}=\frac{1}{\lambda_x}\|\psi_x(\bfx(t),t)\|_{2,[t_0,T]},
\
\gamma_{y_{2,2}}=\frac{1}{\lambda_y}\|\psi_y(\bfy(t),t)\|_{2,[t_0,T]}
\]
respectively, and functions $f_x(\cdot)$, $f_y(\cdot)$ satisfy
Assumptions \ref{assume:alpha}, \ref{assume:alpha_upper} with
\[
\begin{split}
&D_{x}=1.5, \ D_{x,1}=0.5, \ \alpha_x(\bfx,t)= x_1-x_0, \\
&D_{y}=1.6, \ D_{y,1}=0.4, \ \alpha_y(\bfy,t)= y_1-y_0
\end{split}
\]
Hence conditions 1)-4) of Theorem \ref{theorem:interconnection}
are satisfied. Furthermore, according to the remarks regarding
condition 5) of the theorem, requirements
(\ref{eq:disturbance_gain}) can be replaced with implicit
constraints (\ref{eq:disturbance_gain_1}). These, however,
according to (\ref{eq:example_L2_gains}) also hold with
$\beta_x=k_1$, $\beta_y=k_2$.

Given that $\alpha_x(\bfx,t)=x_1-x_0$, $\alpha_y(\bfy,t)=y_1-y_0$,
Assumption \ref{assume:explicit_realizability} will be satisfied
for functions $\alpha_x(\bfx,t)$, $\alpha_y(\bfy,t)$ with
$\Psi_x(\cdot)=0$, $\Psi_y(\cdot)=0$. Therefore, adaptation
algorithms (\ref{fin_forms_ours_tr1x}),
(\ref{fin_forms_ours_tr1y}) will have the following form:
\begin{eqnarray}\label{eq:example_adaptation}
\hat{\theta}_x&=& \Gamma_x((x_1+x_2) (x_1-x_0) +
\hat{\theta}_{x,I}),\nonumber \\
\dot{\hat\theta}_{x,I}&=& \lambda_x (x_1+x_2)(x_1-x_0) - (x_1+x_2)x_2\nonumber \\
\hat{\theta}_y&=& \Gamma_y((y_1+y_2) (y_1-y_0) +
\hat{\theta}_{y,I}),\\
\dot{\hat\theta}_{y,I}&=& \lambda_y (y_1+y_2)(y_1-y_0) -
(y_1+y_2)y_2\nonumber
\end{eqnarray}
Hence, according to Theorem \ref{theorem:interconnection}
boundedness of the solutions in the closed loop system
(\ref{eq:example_error_model}), (\ref{eq:example_adaptation}) is
ensured upon the following condition
\begin{equation}\label{eq:example_condition_boundedness}
\frac{k_1 k_2}{\lambda_x
\lambda_y}\left(1+\frac{D_x}{D_{x,1}}\right)\left(1+\frac{D_y}{D_{y,1}}\right)<1
\Rightarrow k_1 k_2 < \frac{\lambda_x\lambda_y}{20}
\end{equation}
Moreover, given that properties
H\ref{hyp:locally_bound_uniform_f}--
H\ref{hyp:locally_bound_uniform_phi} hold for the chosen functions
$\psi_x(\bfx,t)$, $\psi_y(\bfy,t)$, condition
(\ref{eq:example_condition_boundedness}) guarantees that limiting
relations (\ref{eq:convergence_psi_xy}),
(\ref{eq:convergence_f_xy}) hold.

\begin{figure}
\begin{center}
\includegraphics[width=300pt]{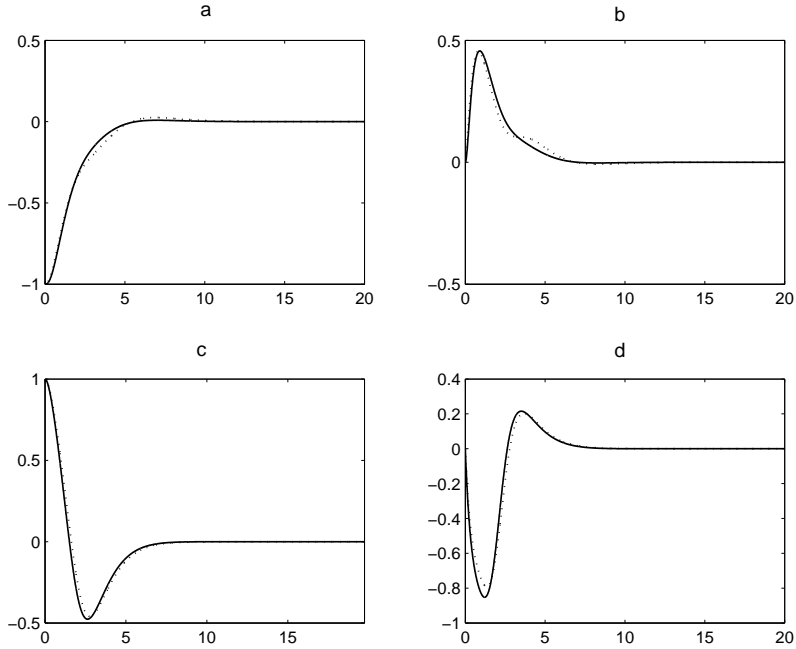}
\end{center}
\begin{center}
\caption{Plots of trajectories $x_1(t)$ (panel a), $x_2(t)$ (panel
b), $y_1(t)$ (panel c), $y_2(t)$ (panel d) as functions of $t$ in
closed loop system (\ref{eq:example_dec_model}),
(\ref{eq:example_control}), (\ref{eq:example_adaptation}). Dotted
lines correspond to the case when $k_1=k_2=0.4$, and solid lines
stand for solutions obtained with the following values of coupling
$k_1=1$, $k_2=0.1$}\label{fig:decentralized:example}
\end{center}
\end{figure}

Trajectories of the closed loop system
(\ref{eq:example_dec_model}), (\ref{eq:example_control}),
(\ref{eq:example_adaptation}) with the following values of
parameters $\Gamma_x=\Gamma_y=1$, $\lambda_x=\lambda_y=2$,
$x_0=y_0=1$, $\theta_x=\theta_y=1$ and initial conditions
$x_1(0)=-1$, $x_2(0)=0$, $y_1(0)=1$, $y_2(0)=0$,
$\hat{\theta}_{x,I}(0)=-1$, $\hat{\theta}_{y,I}(0)=-2$ are
provided in Fig. \ref{fig:decentralized:example}.

\section{Conclusion}

We provided new tools for the  design and analysis of adaptive
decentralized control schemes. Our method allows the desired
dynamics to be Lyapunov unstable and the parametrization of the
uncertainties to be nonlinear. The results are based on a
formulation of the problem for adaptive control as a problem of
regulation in functional spaces (in particular, $L_2^1[t_0,T]$
spaces) rather than of simply reaching of the control goal in
$\Real^n$. This allows us to introduce adaptation algorithms with
new properties and apply a small-gain argument to establish
applicability of these schemes to the problem of decentralized
control.

In order to avoid unnecessary complications, state feedback was
assumed in the main-loop controllers which transform original
equation into the error coupled model. Extension of the results to
output-feedback main loop controllers is a topic for future study.

\section{Proofs of the theorems}

\subsection{Proof of Theorem \ref{stability_theorem}}

Let us first show that property P1) holds. Consider solutions of
system (\ref{system1}), (\ref{error_model_d}),
(\ref{fin_forms_ours_tr1}), (\ref{fin_forms_ours_tr11}) passing
through the point $\bfx(t_0)$, $\hat{\thetavec}_I(t_0)$ for
$t\in[t_0,T^\ast]$
. Let us  calculate the time-derivative  of function
$\hat{\thetavec}(\bfx,t)$:
$\dot{\hat{\thetavec}}(\bfx,t)=\Gamma({\dot{\hat{\thetavec}}_{P}}+\dot{\hat\thetavec}_I)=\Gamma(\dpsi\alphavec(\bfx,t)+\psi\dot{\alphavec}(\bfx,t)-\dot{\Psi}(\bfx,t)+\dot{\hat\thetavec}_I)$.
Notice that
\begin{equation}\label{t2_1}
\begin{split}
&\psi\dot{\alphavec}(\bfx,t)-\dot{\Psi}(\bfx,t)+\dot{\hat{\thetavec}}_I=\psi(\bfx,t)\frac{\pd
\alphavec(\bfx,t)}{\pd \bfx_1}\dot{\bfx}_1+\psi(\bfx,t)\frac{\pd
\alphavec(\bfx)}{\pd \bfx_2}\dot{\bfx}_2 +\\
 & \psi(\bfx,t)\frac{\pd
\alphavec(\bfx,t)}{\pd t}-
 \frac{\pd \Psi(\bfx,t)}{\pd \bfx_1}\dot{\bfx}_1-\frac{\pd
\Psi(\bfx,t)}{\pd \bfx_2}\dot{\bfx}_2-\frac{\pd \Psi(\bfx,t)}{\pd
t}+\dot{\hat\thetavec}_I
\end{split}
\end{equation}
According to Assumption \ref{assume:explicit_realizability},
$\frac{\pd \Psi(\bfx,t)}{\pd \bfx_2}=\psi(\bfx,t)\frac{\pd
\alphavec(\bfx,t)}{\pd \bfx_2}
$. Then taking into account (\ref{t2_1}), we obtain
\begin{equation}\label{t2_2}
\begin{split}
&
\psi\dot{\alphavec}(\bfx,t)-\dot{\Psi}(\bfx,t)+\dot{\hat{\thetavec}}_I=\left(\psi(\bfx,t)\frac{\pd
\alphavec(\bfx,t)}{\pd \bfx_1}-\frac{\pd \Psi}{\pd \bfx_1
}\right)\dot{\bfx}_1\\
&+\psi(\bfx,t)\frac{\pd \alphavec(\bfx,t)}{\pd
t}-\frac{\Psi(\bfx,t)}{\pd t}
\end{split}
\end{equation}
Notice that according to the proposed notation we can rewrite the
term $\left(\psi(\bfx,t)\frac{\pd \alphavec(\bfx,t)}{\pd
\bfx_1}-\frac{\pd \Psi}{\pd \bfx_1 }\right)\dot{\bfx}_1$ in the
following form: $\psi(\bfx,t)L_{\bff_1}
\alphavec(\bfx,t)-L_{\bff_1} \Psi(\bfx,t)+
\left(\psi(\bfx,t)L_{\bfg_1} \alphavec(\bfx,t)-L_{\bfg_1}
\Psi(\bfx,t)\right)u(\bfx,\hat{\thetavec},t)$. Hence, it follows
from (\ref{fin_forms_ours_tr1}) and (\ref{t2_2}) that
$\psi\dot{\alphavec}(\bfx,t)-\dot{\Psi}(\bfx,t)+\dot{\hat{\thetavec}}_I=\varphi(\psi)\alphavec(\bfx,t)
$. Therefore, the derivative $\dot{\hat\thetavec}(\bfx,t)$ can be
written in the following way:
\begin{equation}\label{algorithm_dpsi}
\dot{\hat{\thetavec}}=\Gamma(\dpsi+\varphi(\psi))\alphavec(\bfx,t)
\end{equation}
Asymptotic properties of nonlinear parameterized control systems
with adaptation algorithm (\ref{algorithm_dpsi}) under assumption
of Lyapunov stability of the target dynamics were investigated in
\cite{tpt2003_tac}. In the present contribution we aim to provide
characterizations of the closed loop system in terms of functional
mappings  between functions $\psi(\bfx(t),t)$, $\varepsilon(t)$,
and $f(\bfx(t),\thetavec,t)-f(\bfx(t),\hat{\thetavec}(t),t)$ and
without requiring Lyapunov stability of the target dynamics
(\ref{eq:target_dynamics}).

For this purpose consider the following positive-definite
function:
\begin{equation}\label{V_theta}
V_{\hat{\thetavec}}(\hat{\thetavec},\thetavec,t)=
\frac{1}{2}\|\hat{\thetavec}-\thetavec\|^2_{\Gamma^{-1}} +
\frac{D}{4 D_1^2} \int_{t}^\infty\varepsilon^2(\tau)d\tau
\end{equation}
Its time-derivative according to equations (\ref{algorithm_dpsi})
can be obtained as follows:
\begin{equation}\label{eq:dV_full_alg}
\dot{V}_{\hat{\thetavec}}(\hat{\thetavec},\thetavec,t)=(\varphi(\psi)+\dpsi)(\hat{\thetavec}-\thetavec)^{T}\alphavec(\bfx,t)
-
\frac{D}{4 D_1^2}\varepsilon^2(t)
\end{equation}
%
Hence using Assumptions \ref{assume:alpha},
\ref{assume:alpha_upper} and equality (\ref{error_model_d}) we can
estimate the derivative $\dot{V}_{\hat{\thetavec}}$ as follows:
\begin{eqnarray}\label{parameric_deviation_derivative}
& &
\dot{V}_{\hat{\thetavec}}(\hat{\thetavec},\thetavec,t)\leq-(f(\bfx,\hat{\thetavec},t)-f(\bfx,\thetavec,t)+\varepsilon(t))(\hat{\thetavec}-\thetavec)^{T}\alphavec(\bfx,t)
- \frac{D}{4 D_1^2}\varepsilon^2(t)\nonumber
\\
& &
\leq-\frac{1}{D}(f(\bfx,\hat{\thetavec},t)-f(\bfx,\thetavec,t))^2+\frac{1}{D_1}|\varepsilon(t)||f(\bfx,\hat{\thetavec},t)-f(\bfx,\thetavec,t)|\nonumber\\
& & - \frac{D}{4 D_1^2}\varepsilon^2(t)  \leq -
\frac{1}{D}\left(|f(\bfx,\hat{\thetavec},t)-f(\bfx,\thetavec,t)|-\frac{D}{2
D_1} \varepsilon(t)\right)^2 \leq 0
\end{eqnarray}
It follows immediately from
(\ref{parameric_deviation_derivative}), (\ref{V_theta}) that
\begin{equation}\label{eq:parametric_norm}
\|\hat{\thetavec}(t)-\thetavec\|^{2}_{\Gamma^{-1}}\leq
\|\hat{\thetavec}(t_0)-\thetavec\|^{2}_{\Gamma^{-1}}+\frac{D}{2
D_1^2}\|\varepsilon(t)\|^{2}_{2,[t_0,\infty]}
\end{equation}
In particular, for $t\in[t_0,T^\ast]$ we can derive from
(\ref{V_theta}) that
$\|\hat{\thetavec}(t)-\thetavec\|^{2}_{\Gamma^{-1}}\leq
\|\hat{\thetavec}(t_0)-\thetavec\|^{2}_{\Gamma^{-1}}+\frac{D}{2
D_1^2}\|\varepsilon(t)\|^{2}_{2,[t_0,T^\ast]}$. Therefore
$\hat{\thetavec}(t)\in L_\infty^2[t_0,T^\ast]$. Furthermore
$|f(\bfx(t),\hat{\thetavec}(t),t)-f(\bfx(t),\thetavec,t)|-\frac{D}{2
D_1} \varepsilon(t)\in L_{2}^1 [t_0,T^\ast]$. In particular
\begin{eqnarray}\label{eq:t1_ins1}
&
&\left\||f(\bfx(t),\hat{\thetavec}(t),t)-f(\bfx(t),\thetavec,t)|-\frac{D}{2
D_1} \varepsilon(t)\right\|_{2,[t_0,T^\ast]}^2\leq
\nonumber\\
&&\frac{D}{2}\|\thetavec-\hat{\thetavec}(t_0)\|^{2}_{\Gamma^{-1}}+\frac{D^2}{4
D_1^2}\|\varepsilon(t)\|_{2,[t_0,T^\ast]}^2
\end{eqnarray}
Hence $f(\bfx(t),\hat{\thetavec}(t),t)-f(\bfx(t),\thetavec,t)\in
L_{2}^1 [t_0,T^\ast]$ as a sum of two functions from $L_{2}^1
[t_0,T^\ast]$. In order to estimate the upper bound of the norm
$\|f(\bfx(t),\hat{\thetavec}(t),t)-f(\bfx(t),\thetavec,t)\|_{2,[t_0,T^\ast]}$
from (\ref{eq:t1_ins1}) we use the Minkowski inequality:
\begin{eqnarray}
&&\left\|f(\bfx(t),\hat{\thetavec}(t),t)-f(\bfx(t),\thetavec,t)|-\frac{D}{2
D_1} \varepsilon(t)\right\|_{2,[t_0,T^\ast]}\leq
\nonumber\\
&&\left(\frac{D}{2}\|\thetavec-\hat{\thetavec}(t_0)\|^{2}_{\Gamma^{-1}}\right)^{0.5}+
\frac{D}{2 D_1}\|\varepsilon(t)\|_{2,[t_0,T^\ast]}\nonumber
\end{eqnarray}
and then apply the triangle inequality to the functions from
$L_{2}^1 [t_0,T^\ast]$:
\begin{eqnarray}\label{eq:t1_ins2}
& &
\|f(\bfx(t),\hat{\thetavec}(t),t)-f(\bfx(t),\thetavec,t)\|_{2,[t_0,T^\ast]}\leq\nonumber\\
& &
\left\|f(\bfx(t),\hat{\thetavec}(t),t)-f(\bfx(t),\thetavec,t)-\frac{D}{2
D_1}\varepsilon(t)\right\|_{2,[t_0,T^\ast]}+\\
& & \frac{D}{2 D_1}\|\varepsilon(t)\|_{2,[t_0,T^\ast]}\leq
\left(\frac{D}{2}\|\thetavec-\hat{\thetavec}(t_0)\|^{2}_{\Gamma^{-1}}\right)^{0.5}
+ \frac{D}{D_1}\|\varepsilon(t)\|_{2,[t_0,T^\ast]}\nonumber
\end{eqnarray}
Therefore, property P1) is proven.

Let us prove property P2). In order to do this we have to check
first if the solutions of the closed loop system are defined for
all $t\in\Real_+$, i.e. they do not go to infinity in finite time.
We prove this by a contradiction argument. Indeed, let there
exists time instant $t_s$ such that $\|\bfx(t_s)\|=\infty$. It
follows from P1), however, that
$f(\bfx(t),\hat{\thetavec}(t),t)-f(\bfx(t),\thetavec,t)\in L_{2}^1
[t_0,t_s]$. Furthermore, according to (\ref{eq:t1_ins2}) the norm
$\|f(\bfx(t),\hat{\thetavec}(t),t)-f(\bfx(t),\thetavec,t)\|_{2,[t_0,t_s]}$
can be bounded from above by a continuous function of $\thetavec,
\ \hat{\thetavec}(t_0)$, $\Gamma$,  and
$\|\varepsilon(t)\|_{2,[t_0,\infty]}$. Let us denote this bound by
symbol $D_f$. Notice that $D_f$ does not depend on $t_s$. Consider
system (\ref{error_model_d}) for $t\in[t_0,t_s]$:
$\dpsi=f(\bfx,\thetavec,t)-f(\bfx,\hat{\thetavec},t)-\varphi(\psi,\omegavec,t)+\varepsilon(t)$.
Given that both
$f(\bfx(t),\thetavec,t)-f(\bfx(t),\hat{\thetavec}(t),t),
\varepsilon(t) \in L_{2}^1 [t_0,t_s]$ and taking into account
Assumption \ref{assume:gain}, we automatically obtain that
$\psi(\bfx(t),t)\in L_\infty^{1}[t_0,t_s]$. In particular, using
the triangle inequality and the fact that the function
$\gamma_{\infty,2}\left(\psi(\bfx_0,t_0),\omegavec,M\right)$ in
Assumption \ref{assume:gain} is non-decreasing in $M$, we can
estimate the norm $\|\psi(\bfx(t),t)\|_{\infty,[t_0,t_s]}$ as
follows:
\begin{equation}\label{eq:bound_psi}
\|\psi(\bfx(t),t)\|_{\infty,[t_0,t_s]}\leq
\gamma_{\infty,2}\left(\psi(\bfx_0,t_0),\omegavec,D_f+\|\varepsilon(t)\|^2_{2,[t_0,\infty]}\right)
\end{equation}
According to Assumption \ref{assume:psi} the following inequality
holds:
\begin{equation}\label{eq:bound_x}
\|\bfx(t)\|_{\infty,[t_0,t_s]}\leq\tilde{\gamma}\left(\bfx_0,\thetavec,\gamma_{\infty,2}\left(\psi(\bfx_0,t_0),\omegavec,D_f+\|\varepsilon(t)\|^2_{2,[t_0,\infty]}\right)\right)
\end{equation}
Given that a superposition of locally bounded functions is locally
bounded, we conclude that $\|\bfx(t)\|_{\infty[t_0,t_s]}$ is
bounded. This, however, contradicts to the previous claim that
$\|\bfx(t_s)\|=\infty$. Taking into account inequality
(\ref{eq:parametric_norm}) we can derive that both
$\hat{\thetavec}(\bfx(t),t)$ and $\hat{\thetavec}_I(t)$ are
bounded for every $t\in\Real_+$. Moreover, according to
(\ref{eq:bound_psi}), (\ref{eq:bound_x}),
(\ref{eq:parametric_norm}) these bounds are themselves locally
bounded functions of initial conditions and parameters. Therefore,
$\bfx(t)\in L^n_\infty[t_0,\infty]$,
$\hat{\thetavec}(\bfx(t),t)\in L^d_\infty [t_0,\infty]$.
Inequality (\ref{eq:psi_gain}) follows immediately from
(\ref{eq:t1_ins2}), (\ref{eq:gain_psi_L2}), and the triangle
inequality.  Property P2) is proven.

Let us show that P3) holds. It is assumed that system
(\ref{eq:target_dynamics}) has $L_{2}^1 [t_0,\infty]\mapsto
L_{p}^1 [t_0,\infty]$, $p>1$ gain. In addition, we have just shown
that $f(\bfx(t),\thetavec,t)-f(\bfx(t),\hat{\thetavec}(t),t),
\varepsilon(t) \in L_{2} [t_0,\infty]$. Hence, taking into account
equation (\ref{error_model_d}) we conclude that
$\psi(\bfx(t),t)\in L_{p}^1 [t_0,\infty]$, $p>1$. On the other
hand, given that $f(\bfx,\hat{\thetavec},t)$,
$\varphi(\psi,\omegavec,t)$ are locally  bounded with respect to
their first two arguments uniformly in $t$, and that  $\bfx(t)\in
L_{\infty}^n[t_0,\infty]$,$\psi(\bfx(t),t)\in
L_\infty^1[t_0,\infty]$, $\hat{\thetavec}(t)\in
L_\infty^d[t_0,\infty]$, $\thetavec\in\Omega_\theta$, the signal
$\varphi(\psi(\bfx(t),t),\omegavec,t)+f(\bfx(t),\thetavec,t)-f(\bfx(t),\hat{\thetavec}(t),t)$
is bounded. Then $\varepsilon(t)\in L_\infty^1[t_0,\infty]$
implies that $\dpsi$ is bounded, and P3) is guaranteed by
Barbalat's lemma.

To complete the proof of the theorem  (property P4) consider the
time-derivative of function $f(\bfx,\hat{\thetavec},t)$:
\[
\begin{split}
&\frac{d}{dt}f(\bfx,\hat{\thetavec},t)=L_{\bff(\bfx,\thetavec)+\bfg(\bfx)u(\bfx,\hat{\thetavec},t)}f(\bfx,\hat{\thetavec},t)+\\
& \frac{\pd f(\bfx,\hat{\thetavec},t)}{\pd \hat{\thetavec}}\Gamma
(\varphi(\psi,\omegavec,t)+\dpsi)\alphavec(\bfx,t)+\frac{\pd
f(\bfx,\hat{\thetavec},t)}{\pd t}
\end{split}
\]
Taking into account that the function $f(\bfx,\thetavec,t)$ is
continuously differentiable in $\bfx$, $\thetavec$; the derivative
$ \pd {f(\bfx,\thetavec,t)}/{\pd t}$ is locally bounded with
respect to $\bfx$, $\thetavec$ uniformly in $t$; functions
$\alphavec(\bfx,t)$, $\pd \psi(\bfx,t)/\pd t$ are locally bounded
with respect to $\bfx$ uniformly in $t$, then $d/dt
(f(\bfx,\thetavec,t)-f(\bfx,\hat\thetavec,t))$ is bounded. Then
given that
$f(\bfx(t),\thetavec,t)-f(\bfx(t),\hat{\thetavec}(t),t)\in L_{2}^1
[t_0,\infty]$ by applying Barbalat's lemma we conclude that
$f(\bfx,\thetavec,\tau)-f(\bfx,\hat\thetavec,\tau)\rightarrow 0$
as $t\rightarrow\infty$. { The theorem is proven.}

\subsection{Proof of Theorem \ref{theorem:interconnection}}

Let us denote
\[
\Delta f_x[t_0,T]=
\|f_x(\bfx,\thetavec_x,t)-f_x(\bfx,\hat{\thetavec}_x,t)\|_{2,[t_0,T]},
\]
\[
\Delta f_y
[t_0,T]=\|f_x(\bfy,\thetavec_y,t)-f_y(\bfy,\hat{\thetavec}_y,t)\|_{2,[t_0,T]}.
\]
As follows from Theorem \ref{stability_theorem} the following
inequalities hold
\begin{equation}\label{proof:interconnection:t1}
\Delta f_x[t_0,T]\leq C_x + \frac{D_x}{D_{1,x}}
\|h_y(\bfx(t),\bfy(t),t)\|_{2,[t_0,T]}
\end{equation}
\begin{equation}\label{proof:interconnection:t2}
\Delta f_y[t_0,T]\leq C_y + \frac{D_y}{D_{1,y}}
\|h_x(\bfx(t),\bfy(t),t)\|_{2,[t_0,T]},
\end{equation}
where $C_x$, $C_y$ are some constants, independent of $T$. Taking
estimates (\ref{proof:interconnection:t1}),
(\ref{proof:interconnection:t2}) into account  we obtain the
following estimates:
\begin{equation}\label{proof:interconnection:t3}
\begin{split}
&\Delta f_x[t_0,T]+\|h_y(\bfx(t),\bfy(t),t)\|_{2,[t_0,T]}\leq \\
&C_x + \left(\frac{D_x}{D_{1,x}}+1\right)
\|h_y(\bfx(t),\bfy(t),t)\|_{2,[t_0,T]}
\end{split}
\end{equation}
\begin{equation}\label{proof:interconnection:t4}
\begin{split}
&\Delta f_y[t_0,T]+\|h_x(\bfx(t),\bfy(t),t)\|_{2,[t_0,T]}\leq
\\
&C_y + \left(\frac{D_y}{D_{1,y}}+1\right)
\|h_x(\bfx(t),\bfy(t),t)\|_{2,[t_0,T]},
\end{split}
\end{equation}
The proof of the theorem would be complete if we show that the
$L_2^1[t_0,T]$ norms of $h_x(\bfx(t),\bfy(t),t)$,
$h_y(\bfx(t),\bfy(t),t)$ are globally bounded uniformly in $T$.
Let us show that this is indeed the case. Using the widely known
generalized triangular inequality \cite{Jiang_1994}
\[
\gamma(a + b)\leq \gamma((\rho+Id)(a))+\gamma((\rho+Id)\circ
\rho^{-1}(b)), \ a,b\in\Real_+, \ \gamma,\rho\in\mathcal{K}_\infty,
\]
equations (\ref{proof:interconnection:t3}),
(\ref{proof:interconnection:t4}) and also property
(\ref{eq:disturbance_gain}), we conclude that
\begin{equation}\label{proof:interconnection:t5}
\begin{split}
&\|h_y(\bfx(t),\bfy(t),t)\|_{2,[t_0,T]}\leq\\
&\beta_y\cdot \gamma_{y_{2,2}}\circ\rho_1
\left(\left(\frac{D_y}{D_{1,y}}+1\right)\|h_x(\bfx(t),\bfy(t),t)\|_{2,[t_0,T]}\right)+C_{y,1}\\
&\|h_x(\bfx(t),\bfy(t),t)\|_{2,[t_0,T]}\leq\\
& \beta_x\cdot \gamma_{x_{2,2}}\circ\rho_2
\left(\left(\frac{D_x}{D_{1,x}}+1\right)\|h_y(\bfx(t),\bfy(t),t)\|_{2,[t_0,T]}\right)+C_{x,1}
\end{split}
\end{equation}
where $\rho_1(\cdot)$, $\rho_2(\cdot)\in\mathcal{K}_\infty$,
$\rho_1(\cdot), \rho_2(\cdot)>Id(\cdot)$. Then, according to
(\ref{proof:interconnection:t5}), the existence of
$\rho_3(\cdot)\in\mathcal{K}_\infty\geq Id(\cdot)$, satisfying
inequality
\[
\beta_y\circ\gamma_{y_{2,2}}\circ\rho_1\circ\left(\frac{D_y}{D_{y,1}}+1\right)\circ\rho_3\circ
\beta_x\circ
\gamma_{x_{2,2}}\circ\rho_2\circ\left(\frac{D_x}{D_{x,1}}+1\right)(\Delta)<
\Delta \ \forall \ \Delta\geq \bar{\Delta}
\]
for some $\bar{\Delta}\in\Real_+$ ensures that the norms
\[
\|h_y(\bfx(t),\bfy(t),t)\|_{2,[t_0,T]}, \
\|h_x(\bfx(t),\bfy(t),t)\|_{2,[t_0,T]}
\]
are globally uniformly bounded in $T$.  The rest of the proof
follows from Theorem \ref{stability_theorem}.  { The theorem is
proven.}

\bibliographystyle{plain}
\bibliography{Decentralized}

\end{document}